\documentclass[twoside,11pt]{article}

\pdfoutput=1   

\usepackage[utf8]{inputenc}   
\usepackage[T1]{fontenc}      
\usepackage[english]{babel}
\usepackage{amsfonts} 
\usepackage{amsmath}
\usepackage{latexsym} 
\usepackage{lmodern}  
\usepackage[all]{xy}
\usepackage{color}
\usepackage{mathrsfs}\let\cal\mathscr

\title{A Construction of Colimits\\ in\\ Monoidal Closed Categories}
\author{Alain Prouté\thanks{Université Paris-Diderot - France}}   

\SilentMatrices
\nonstopmode   
\DeclareMathOperator{\Lim}{{\bf lim}}
\DeclareMathOperator{\ev}{{{\bf\,ev}}}

\DeclareMathOperator{\op}{{\rm op}}
\DeclareMathOperator{\un}{{\bf 1}}

\DeclareMathOperator{\Ab}{{\mbox{\bf Ab}}}
\newcommand \la {\langle}
\newcommand \ra {\rangle}
\newcommand \Z {{\mathbb Z}}
\newcommand \CC {{\cal C}}
\newcommand \DD {{\cal D}}
\newcommand \II {{\cal I}}
\newcommand \rar {{\xymatrix@1@C=4mm{{}\ar[r]&{}}}}
\newcommand \ci {\circ}
\newcommand \ba[1] {\overline{#1}}
\newcommand \qed {\mbox{$\Box$}}
\newtheorem{theorem}{Theorem}{}
\newtheorem{lemma}{Lemma}{}

\begin{document}
\everymath{\displaystyle}

\maketitle

   \begin{abstract}
   We prove that a
   category which is symmetric (relaxed) monoidal closed, (small) complete, 
   well-powered and has a small cogenerating family, is (small) cocomplete. 
   \end{abstract}

\section{Monoidal closed categories}
For the notion of (relaxed) symmetric monoidal closed category, I refer to the
literature, for example Mac~Lane \cite{MacLane63} and \cite{MacLane}. 
In this section, I recall some
facts about this structure mainly for fixing the notations. 

A monoidal category is equiped with a functor $\otimes:\CC\times \CC\rar \CC$, and a
particular object $I$, together with isomorphisms:
$$
\xymatrix@R=0mm@C=14mm{
X\otimes (Y\otimes Z)\ar[r]^-{\alpha_{X,Y,Z}}&(X\otimes Y)\otimes Z\\
X\otimes I\ar[r]^-{\rho_X}&X\\
I\otimes X\ar[r]^-{\lambda_X}&X
}
$$
natural in $X$, $Y$ and $Z$, and satisfying some axioms \cite{MacLane}. 
Such a category is said to be \emph{symmetric} if there is an
isomorphism
$$
\xymatrix{
X\otimes Y\ar[r]^-{\sigma_{X,Y}}&Y\otimes X
}
$$
natural in $X$ and $Y$, satisfying the equation $\sigma_{Y,X}\ci\sigma_{X,Y} =
1_{X\otimes Y}$, $\lambda_X\ci\sigma_{X,I}=\rho_X$, 
and some more axioms. 

A symmetric monoidal category is \emph{closed} if for each object $Y$, the
functor $X\mapsto X\otimes Y$ has a right adjoint, which is denoted $Z\mapsto
Z^Y$. In other words there is a bijection 
$$
\xymatrix@C=10mm{
\CC(X\otimes Y,Z)\ar[r]_-\simeq^-{\Lambda_Y}&\CC(X,Z^Y)
}
$$
which is natural in $X$ and $Z$. 
The co-unit of the adjunction is an arrow
$\ev:Z^Y\otimes Y \rar Z$.

The correspondance $(Y,Z)\mapsto Z^Y$ is actually a bifunctor, i.e. a functor
$\CC^{\op}\times\CC\rar\CC$, so that for any arrow $f:Y\rar Y'$, we have an
arrow $Z^f:Z^{Y'}\rar Z^Y$, and for any arrow $g:Z\rar Z'$, we have an arrow
$g^Y:Z^Y\rar {Z'}^Y$. 

An arrow $e:I\rar X$ will be called an \emph{element} of $X$.(\footnote{Notice
that $I$ is not necessarily a terminal object, for example in the category $\Ab$
of abelian groups with usual tensor product, where $I$ is $\Z$ and where the
terminal object is $0$.}) For any object $X$, we have a canonical isomorphism
$\iota_X = \Lambda_X(\rho_X):X\rar X^I$,  
and a canonical arrow $\eta_X = \Lambda_X(\lambda_X):I\rar X^X$.

To any arrow $f:X\rar Z^Y$ corresponds the arrow $\ba f:Y\rar Z^X$, defined by 
$\ba f = \Lambda_X(\ev\ci(f\otimes 1_Y)\ci\sigma_{Y,X})$. We call $\ba f$ the
\emph{swap} of $f$. It is easily checked that $\ba{\ba f} = f$. The swap also has the
following easily checked properties:
$$
\begin{array}{rclll}
\iota_X &=& \ba{\eta_X}\\
Z^f\ci \ba g &=& \ba{g\ci f} &\hspace{8mm}&\mbox{for
    $\xymatrix@1{X\ar[r]^-f&X'\ar[r]^-g&Z^Y}$}\\
g^X\ci\ba f &=& \ba{g^Y\ci f}&&\mbox{for  $\xymatrix@1{X\ar[r]^-f&Z^Y}$ and 
                                          $\xymatrix@1{Z\ar[r]^-g&Z'}$}
\end{array}
$$

If we have an element $e:I\rar X$, then for any object $Y$, we have an
\emph{evaluation at $e$} arrow $\ev_e = \iota_Y^{-1}\ci Y^e:Y^X\rar Y$.
$$
\xymatrix{
Y^X\ar@(dr,dl)[rr]_-{\ev_e}\ar[r]^-{Y^e}&Y^I\ar[r]^-{\iota_Y^{-1}}&Y
}
$$

\section{Limits, colimits and ends}

The notions of \emph{limit}, \emph{colimit} and \emph{end} are assumed to be 
known by the reader. Full details can be found in
Mac~Lane \cite{MacLane}. Nevertheless, we recall some basic facts because we 
need to fix some notations.

Given a diagram $d:\II\rar\CC$, a \emph{cone} over $d$ with vertex $X$, is a
family of arrows $(\rho_i:X\rar d(i))_i$ (called \emph{edges} or \emph{canonical 
    projections}), one for each object in $\II$, 
such that for each arrow $\theta:i\rar j$ of $\II$, the diagram:
$$
\xymatrix@R=2mm{
d(i)\ar[dd]_-{d(\theta)}\\
&X\ar[ul]_-{\rho_i}\ar[dl]^-{\rho_j}\\
d(j)
}
$$
is commutative. Given two cones on $d$, $(\rho_i:X\rar d(i))_i$ and $(\pi_i:Y\rar d(i))_i$, a morphism
from the former to the later is an arrow $f:X\rar Y$ such that 
$$
\xymatrix@R=4mm@C=12mm{
d(i)\ar[dd]_-{d(\theta)}\\
&X\ar[ul]^-{\rho_i}\ar[dl]_-{\rho_j}\ar[r]^(.3)f&Y\ar@/_/[ull]_-{\pi_i}\ar@/^/[dll]^-{\pi_j}\\
d(j)
}
$$
is commutative. A \emph{limiting cone} over $d$ is a cone which is terminal among
all cones over $d$. Hence, assuming that $(\pi_i:Y\rar d(i))_i$ is a limiting cone in the
above diagram, we have one and only arrow $f:X\rar Y$ such that $\pi_i\ci f =
\rho_i$ for all $i$. This arrow will be
denoted $\la \rho_i\ra_i$. We have $\pi_i\ci\la \rho_i\ra_i = \rho_i$ and 
$\la \rho_i\ra_i\ci\varphi = \la \rho_i\ci\varphi\ra_i$, for any $\varphi$. The
object $Y$ itself is denoted $\Lim d$. 

Actually, $\Lim:\CC^\II\rar\CC$ is a functor (which is right adjoint to the
    diagonal functor $\Delta:\CC\rar\CC^\II$). Hence for any natural
transformation $\varphi:d\rar d'$ between two $\II$-diagrams in $\CC$, we have an arrow
$\Lim \varphi:\Lim d\rar\Lim d'$.

The notion of \emph{colimit} is dual to the notion of limit, i.e. there is a
notion of \emph{cocone} (same as a cone but with the edges in the opposite
direction), and a \emph{colimiting cocone} of $d$ is initial among all cocones
over $d$.

Any cocone $\delta = (\delta_i:d(i)\rar X)_i$ over $d$, gives rise to an element in $\Lim
X^d$. Indeed, we have the cone $(X^{\delta_i}:X^X\rar
X^{d(i)})_i$ over $X^d$, hence the arrow $\la X^{\delta_i}\ra_i: X^X\rar \Lim
X^d$, so that the composition
$$
\xymatrix@C=12mm{
I\ar[r]^-{\eta_X}&X^X\ar[r]^-{\la X^{\delta_i}\ra_i}&\Lim X^d
}
$$
is the wanted element, that we denote $[\delta]$. Because of this, we have an
evaluation arrow $\ev_{[\delta]}:X^{\Lim X^d}\rar X$, and we have:
$
\ev_{[\delta]}\ci\ba{\pi_{i,X}} = \delta_i
$
where $\pi_{i,X}:\Lim X^d\rar X^{d(i)}$ is the canonical projection. 
Indeed:
$$
\begin{array}{rcl}
\ev_{[\delta]}\ci\ba{\pi_{i,X}} 
   &=& \iota_X^{-1}\ci X^{[\delta]}\ci\ba{\pi_{i,X}}\\
   &=& \iota_X^{-1}\ci\ba{\pi_{i,X}\ci[\delta]}\\
   &=& \iota_X^{-1}\ci\ba{\pi_{i,X}\ci \la X^{\delta_i}\ra_i\ci\eta_X}\\
   &=& \iota_X^{-1}\ci \ba{X^{\delta_i} \ci \eta_X}\\
   &=& \iota_X^{-1}\ci \ba{\eta_X} \ci \delta_i\\
   &=& \delta_i
\end{array}
$$

An \emph{end} is a special kind of limit. A \emph{bifunctor} from $\CC$ to $\DD$
is a functor $B:\CC^{\op}\times\CC\rar \DD$, i.e. a functor of \emph{two
variables}, which is contravariant in the first one and covariant in the second
one. An \emph{end} of $B$, denoted $\int_X B(X,X)$, is the vertex of a limiting cone over the
\emph{subdivision diagram} associated to $B$. This diagram has one object $B(X,X)$ for
each object $X$ and one object $B(X,Y)$ for each arrow $f:X\rar Y$ of $\CC$, and 
the two arrows shown below for each arrow $f:X\rar Y$ of $\CC$:
$$
\xymatrix{
B(X,X)\ar[dr]_-{B(X,f)}&&B(Y,Y)\ar[dl]^-{B(f,Y)}\\
&B(X,Y)
}
$$
Consequently, for each object $X$ of $\CC$, we have the edge $\pi_X:\left(\int_X
B(X,X)\right)\rar B(X,X)$ and the relation $B(X,f)\ci\pi_X = B(f,Y)\ci\pi_Y$, the two
members of this equation being actually the edge targeting the instance of 
$B(X,Y)$ associated to $f$. 

Given an object $Z$ of $\CC$ and a family of arrows $g_X:Z\rar B(X,X)$, one for
each object $X$ of $\CC$, such that $B(X,f)\ci g_X = B(f,Y)\ci g_Y$ for each
$f:X\rar Y$ in $\CC$, there is one and only one arrow $\varphi:Z\rar \int_X
B(X,X)$ such that $\pi_X\ci\varphi = g_X$ for all $X$. This arrow is denoted
$\la g_X\ra_X$, where the last $X$ mutes the first one. So, we have the 
relations $\pi_X\ci\la g_X\ra_X = g_X$, and also $\la g_X\ra_X\ci\varphi = \la
g_X\ci\varphi\ra_X$ for any $\varphi$.

\section{Equivalent diagrams}

We say that two diagrams $d_1,d_2;\II\rar \CC$ are \emph{isomorphic} if the
two functors $d_1$ and $d_2$ are isomorphic, i.e. if there exists an isomorphism
$d_1(i)\rar d_2(i)$ natural in $i$. 

\emph{Remark:} It is immediate that if $d_1$ is
isomorphic to $d_2$ and if $d_1$ has a limit, then $d_2$ also has a limit
with the same vertex. 

We say that two diagrams $d_1:\II_1\rar\CC$ and $d_2:\II_2\rar\CC$ are
\emph{equivalent} if there is an equivalence of categories $\Phi:\II_1\rar
\II_2$ such that the two diagrams $d_2\ci\Phi$ and $d_1$ are isomorphic.  
$$
\xymatrix{
\II_1\ar[d]^-{\Phi}_-\simeq\ar[r]^-{d_1}&\CC\\
\II_2\ar[ru]_{d_2}
}
$$
\emph{Remark:} Let $\gamma_i:d_2(\Phi(i))\rar d_1(i)$ be an isomorphism natural
in $i$, and let
$\Psi$ be an inverse of $\Phi$ up to natural isomorphism, so that we have an
isomorphism $\alpha_j:\Phi(\Psi(j))\rar j$, natural in $j$.
Then, the composition
$$
\xymatrix@C=14mm{
d_1(\Psi(j))\ar[r]^-{\gamma_{\Psi(j)}^{-1}}&d_2(\Phi(\Psi(j)))\ar[r]^-{d_2(\alpha_j)}&d_2(j)
}
$$
is an
isomorphism $\delta_j:d_1(\Psi(j))\rar d_2(j)$ which is natural in $j$.

\begin{lemma}
If the diagrams $d_1:\II_1\rar\CC$ and $d_2:\II_2\rar\CC$ are equivalent, and if
$d_1$ has a limit, then so does $d_2$ with the same vertex. 
\end{lemma}

Because of the first remark above, we can assume that $d_1=d_2\ci\Phi$.
Let $(\mu_i:L\rar d_1(i))_i$ be a limiting cone on $d_1$. 
Let $\lambda:j\rar k$ be an arrow of $\II_2$. We have the
commutative diagram:
$$
\xymatrix{
&L\ar[dl]_-{\mu_{\Psi(j)}}\ar[dr]^-{\mu_{\Psi(k)}}\\
d_1(\Psi(j))\ar[d]_-\simeq^-{{\delta_{j}}}
              \ar[rr]_-{d_1(\Psi(\lambda))}&&d_1(\Psi(k))\ar[d]^-\simeq_-{{\delta_{k}}}\\
d_2(j)\ar[rr]_-{d_2(\lambda)}&&d_2(k)
}
$$
We want to prove that $({\delta_j}\ci\mu_{\Psi(j)}:L\rar d_2(j))_j$ is a limiting
cone on $d_2$. Let $(\nu_j:X\rar d_2(j))_j$ be an arbitrary cone on $d_2$, so
that we have $d_2(\lambda)\ci\nu_j = \nu_k$ for all arrows $\lambda:j\rar k$ in
$\II_2$, in particular for the arrow $\alpha_k:\Phi(\Psi(k))\rar k$. 

For any object
$i$ in $\II_1$, we have the composition:
$$
\xymatrix@C=12mm{
X\ar[r]^-{\nu_{\Phi(i)}}&d_2(\Phi(i))\ar[r]^-{\gamma_i}&d_1(i)
}
$$
Because $\gamma_i$ is natural in $i$, these arrows make a cone on $d_1$, so that we have a
unique arrow $f:X\rar L$ such that $\mu_i\ci f = \gamma_i\ci\nu_{\Phi(i)}$ for
any object $i$ of $\II_1$. Hence, for any object $k$ of $\II_2$, we have:
$$
\begin{array}{rcl}
{\delta_k}\ci\mu_{\Psi(k)}\ci f 
  &=&d_2(\alpha_k)\ci\gamma_{\Psi(k)}^{-1}\ci\gamma_{\Psi(k)}\ci\nu_{\Phi(\Psi(k))}\\
  &=&d_2(\alpha_k)\ci\nu_{\Phi(\Psi(k))}\\
  &=&\nu_k
\end{array}
$$
It remains to prove the uniqueness of $f$ such that
${\delta_k}\ci\mu_{\Psi(k)}\ci f = \nu_k$. If for all $k$ we have 
${\delta_k}\ci\mu_{\Psi(k)}\ci f = {\delta_k}\ci\mu_{\Psi(k)}\ci g$, we also
have $\mu_{\Psi(k)}\ci f = \mu_{\Psi(k)}\ci g$ because $\delta_k$ is an
isomorphism. So that for any object $i$ of $\II_1$, we have
$\mu_{\Psi(\Phi(i))}\ci f = \mu_{\Psi(\Phi(i))}\ci g$. But since we also have an
isomorphism $\beta_i:\Psi(\Phi(i))\rar i$ natural in $i$, we have 
$d_1(\beta_i)\ci\mu_{\Psi(\Phi(i))}\ci f = d_1(\beta_i)\ci\mu_{\Psi(\Phi(i))}\ci
g$. i.e. $\mu_i\ci f = \mu_i\ci g$ for all $i$, so that $f=g$.~\qed

\section{On the existence of ends}

\emph{Remark:} Before proving the lemma below, I feel necessary to recall some facts about how
diagrams are graphically represented. By definition, a diagram in $\CC$ is a
(covariant) functor $d:\II\rar\CC$. When we draw diagrams, we don't represent 
the functor $d$ itself but only
its image in $\CC$. As an example, consider the
categories (where identity arrows are not represented):
$$
\II_1 = \xymatrix{i\ar@/^/[r]^-f\ar@/_/[r]_-g&j}
\hspace{10mm}
\II_2 = \raisebox{17pt}{\xymatrix@R=0mm{&k\\
                  i\ar[ur]^-f\ar[dr]_-g\\
                  &l}}
$$
The limit of an $\II_1$-diagram $d$ is an equalizer of $d(f)$ and $d(g)$.
However, a limit of an $\II_2$-diagram $d$ is an isomorphism $\bullet\rar d(i)$,
\emph{even if $d(k)=d(l)$}.
In order to avoid any confusion that could be fatal
for the correctness of a proof, it is important that for any distinct objects
$i$ and $j$ of $\II$, the objects $d(i)$ and $d(j)$ have separate
representations. One can be tempted to represent them at the same position in
the diagram, but clearly this can lead to a misinterpretation. 

\begin{lemma}
Let $F:\CC\rar\CC$ be a (covariant) endofunctor of a symmetric monoidal closed complete 
well-powered category $\CC$ with a small cogenerating family. Then the end
$$
\int_X X^{F(X)}
$$
exists. 
\end{lemma}
Since $\CC$ is complete, the product $P = \Pi G$ of all objects of the cogenerating
family exists and is a cogenerator. For each object $X$ of
$\CC$, we consider the limiting cone with vertex $M_X$:
$$
\xymatrix{
                   &P^{F(P)}\ar[dr]_-{P^{F(\varphi)}}\ar[drrr]^-{P^{F(\psi)}}\\
M_X\ar[ur]^-{m_X}\ar[dr]_-{n_X}  &           &P^{F(X)} &\dots   &P^{F(X)}\\
                   &X^{F(X)}\ar[ur]^-{\varphi^{F(X)}}\ar[urrr]_-{\psi^{F(X)}}
}
$$
where all arrows $\varphi,\psi,\dots:X\rar P$ intervene. We claim that the arrow
$m_X:M_X\rar P^{F(P)}$ is a monomorphism. Indeed, let $u,v:Z\rar M_X$ be two
arrows such that $m_X\ci u = m_X\ci v$. We have $\varphi^{F(X)}\ci n_X\ci u =
\varphi^{F(X)}\ci n_X\ci v$ for all arrows $\varphi:X\rar P$.  Because $Y\mapsto Y^{F(X)}$
is a right adjoint, it preserves monomorphic families. Hence, we have 
$n_X\ci u = n_X\ci v$, which together with $m_X\ci u = m_X\ci v$ entails $u =
v$. 

We consider the diagram made of all arrows $m_X$ and all isomorphisms
$r_{X,Y}:M_X\rar M_Y$ such that $m_Y\ci r_{X,Y} = m_X$. Because $\CC$ is
well-powered, this big diagram is equivalent to a small one, so that it has a
limiting cone $(\zeta_X:\Gamma\rar M_X)_X$. Then, 
$(\pi_X = n_X\ci\zeta_X:\Gamma\rar X^{F(X)})_X$ is a cone on the
subdivision diagram associated to the bifunctor $(X,Y)\mapsto Y^{F(X)}$.  
Indeed, consider the diagram below, where $f:X\rar Y$:
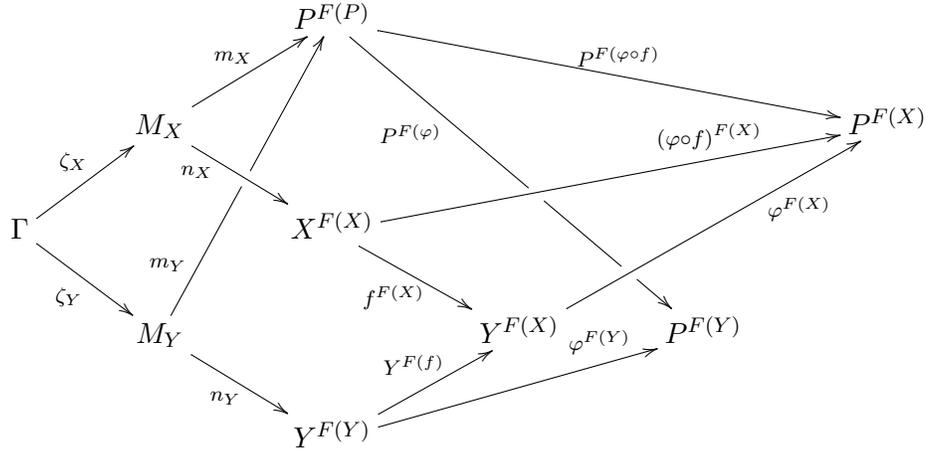
$$
\xymatrix@R=8mm@C=12mm{
                   &&P^{F(P)}\ar[drrr]^-{P^{F(\varphi\ci f)}}
                            \ar[rrddd]|(.55)\hole|(.8)\hole_(.3){P^{F(\varphi)}}\\
&M_X\ar[ur]^-{m_X}\ar[dr]_(.3){n_X}  &      & &  &P^{F(X)}\\
\Gamma\ar[ur]^{\zeta_X}\ar[dr]_-{\zeta_Y}    &&X^{F(X)}\ar[dr]_-{f^{F(X)}}\ar[urrr]^(.70){{(\varphi\ci f)}^{F(X)}}\\
&M_Y\ar[uuur]^(.2){m_Y}|\hole\ar[dr]_-{n_Y}
      &&Y^{F(X)}\ar[rruu]_(.7){\varphi^{F(X)}}&P^{F(Y)}\\
&&Y^{F(Y)}\ar[ur]^-{Y^{F(f)}}\ar[urr]^(.75){\varphi^{F(Y)}}
}
$$
We must prove that the bottom left hexagon is commutative. By construction, we
have $\varphi^{F(X)}\ci Y^{F(f)}\ci n_Y\ci\zeta_Y = \varphi^{F(X)}\ci
f^{F(X)}\ci n_X\ci \zeta_X$, for all arrows $\varphi:Y\rar P$. As previously,
the set of arrows $\varphi^{F(X)}:Y^{F(X)}\rar P^{F(X)}$ is a monomorphic
family, so that our diagram is commutative.

It remains to prove that the cone $(\pi_X:\Gamma\rar X^{F(X)})_X$ is terminal
among the cones over the subdivision diagram. Let $(\xi_X:U\rar X^{F(X)})_X$ be
an arbitrary cone over this diagram. Because $P^{F(\varphi)}\ci \xi_P =
\varphi^{F(Y)}\ci \xi_Y$ for all $\varphi:Y\rar P$, the arrow $\xi_Y$ has a
lifting along $n_Y$ as an arrow $\theta_Y:U\rar M_Y$, and we have
$m_Y\ci\theta_Y = \xi_P$.

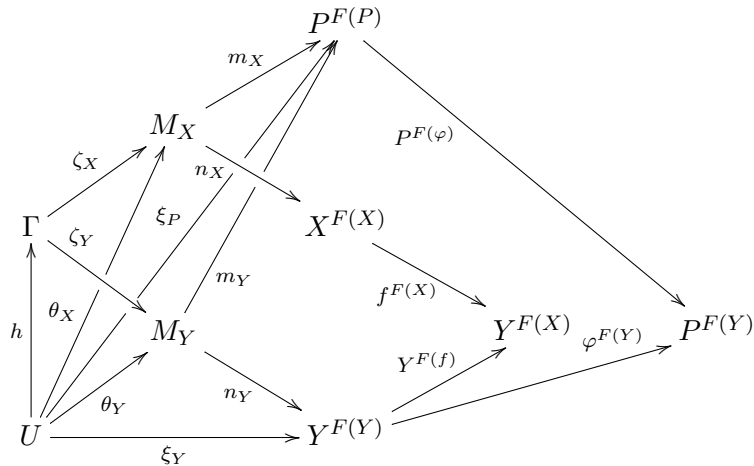
$$
\xymatrix@R=8mm@C=12mm{
                   &&P^{F(P)}\ar[rrddd]_(.3){P^{F(\varphi)}}\\
&M_X\ar[ur]^-{m_X}\ar[dr]_(.3){n_X}\\
\Gamma\ar[ur]^{\zeta_X}\ar[dr]^(.25){\zeta_Y}    
         &&X^{F(X)}\ar[dr]_-{f^{F(X)}}\\
&M_Y\ar[uuur]_(.2){m_Y}|\hole\ar[dr]_-{n_Y}
      &&Y^{F(X)}&P^{F(Y)}\\
U\ar[ur]_-{\theta_Y}\ar[uuur]|\hole^(.38){\theta_X}\ar[rr]_-{\xi_Y}\ar[uu]^-h
 \ar[rruuuu]|(.32)\hole^-{\xi_P}|(.66)\hole
      &&Y^{F(Y)}\ar[ur]^-{Y^{F(f)}}\ar[urr]^(.75){\varphi^{F(Y)}}
}
$$
We have $m_Y\ci\theta_Y = \xi_P = m_X\ci\theta_X$ for all $X$ and $Y$, and if 
$r_{X,Y}:M_X\rar M_Y$ is an isomorphism such that $m_Y\ci r_{X,Y} = m_X$, we
have $r_{X,Y}\ci\theta_X = \theta_Y$ since $m_Y$ is a monomorphism. 

Consequently, there is a unique arrow $h:U\rar\Gamma$ such that $\zeta_Y\ci h =
\theta_Y$ for all $Y$, so that we have 
$\pi_Y\ci h = n_Y\ci\zeta_Y\ci h = n_Y\ci\theta_Y = \xi_Y$. 

It remains to prove the uniqueness of $h$ such that $\pi_Y\ci h = \xi_Y$ for all
$Y$. Let $k:U\rar \Gamma$ be such that $\pi_Y\ci k = \xi_Y$ for all $Y$. We have
$\varphi^{F(Y)}\ci n_Y\ci \zeta_Y\ci k = \varphi^{F(Y)}\ci\xi_Y =
P^{F(\varphi)}\ci m_Y\ci\theta_Y$, so that by definition of $M_Y$, we have 
$\zeta_Y\ci k = \theta_Y$, which proves that $k=h$.~\qed

\section{Statement and proof of the theorem}

\begin{theorem}\label{thethm}
Let $\CC$ be a (locally small) symmetric (relaxed) monoidal closed well-powered category, 
which is (small) complete and 
   has a small cogenerating family. Then, $\CC$ is (small) cocomplete.
\end{theorem}

\emph{Proof of the theorem.} Let $d:\II\rar\CC$ be a small diagram in $\CC$
(i.e. $\II$ is a small category). We have to construct an initial object in the
category of cocones over $d$. This category is actually a comma-category,
precisely the category $\ba{d}/\Delta$, where $\ba{d}:\un\rar \CC^\II$ is the
swap of $d:\II\rar\CC^{\un}$ (after identification of $\CC$ with $\CC^{\un}$),
and where $\Delta:\CC\rar\CC^\II$ is the diagonal functor. Since $\CC$ is
complete, so is $\ba{d}/\Delta$, which is also locally small. Hence, by Freyd's
initial object theorem (see \cite{MacLane}), it is enough to prove that there is
a weakly initial cocone on $d$, i.e. a cocone from which there is at least one
morphism to any other cocone on $d$, but not necessarily at most one. It is
furthermore immediate that the unique morphism from an initial object towards a
weakly initial object has a retraction.

Given an object $X$ of $\CC$, we can consider the diagram $X^d$
in $\CC$, which is actually a functor $X^d:\II^{\op}\rar\CC$, mapping $i$ to
$X^{d(i)}$ and $\theta:i\rar j$ to $X^{d(\theta)};X^{d(j)}\rar X^{d(i)}$. Since
$\CC$ is small-complete, we can consider the objects $\Lim X^d$ and 
$X^{\Lim X^d}$ of $\CC$. 

By the lemma above, we can define:
$$
\Gamma = \int_X X^{\Lim X^d}
$$
We shall prove that $\Gamma$ is the vertex of a weakly initial cocone on $d$. 

Our first task is to define the edges of the cocone with vertex $\Gamma$, that must be arrows
$\gamma_i:d(i)\rar\Gamma$. We have the arrow $\pi_{i,X}:\Lim X^d\rar X^{d(i)}$
which is a canonical projection
of the limiting cone, whose swap is $\ba{\pi_{i,X}}:d(i)\rar X^{\Lim X^d}$.
Furthermore, for any $f:X\rar Y$, we have the commutative diagram:
$$
\xymatrix@R=2mm@C=12mm{
&X^{\Lim X^d}\ar[dr]^-{f^{\Lim X^d}}\\
d(i)\ar[ur]^-{\ba{\pi_{i,X}}}\ar[dr]_-{\ba{\pi_{i,Y}}}&&Y^{\Lim X^d}\\
&Y^{\Lim Y^d}\ar[ur]_-{Y^{\Lim f^d}}
}
$$
Indeed, we have $f^{\Lim X^d}\ci\ba{\pi_{i,X}} = \ba{f^{d(i)}\ci\pi_{i,X}}$ on
the one hand, and $Y^{\Lim f^d}\ci\ba{\pi_{i,Y}} = \ba{\pi_{i,Y}\ci\Lim f^d} =
\ba{f^{d(i)}\ci\pi_{i,X}}$ on the other hand. This ensures that we have the
arrow $\gamma_i = \la \ba{\pi_{i,X}}\ra_X:d(i)\rar\Gamma$.

Now, we must check that for any $\theta:i\rar j$, we have $\gamma_j\ci d(\theta)
= \gamma_i$. For any $X$, we have the commutative diagram:
$$
\xymatrix@R=2mm@C=12mm{
d(i)\ar[dd]_-{d(\theta)}\ar[dr]^-{\ba{\pi_{i,X}}}\\
&X^{\Lim X^d}\\
d(j)\ar[ur]_-{\ba{\pi_{j,X}}}
}
$$
Indeed, $\ba{\pi_{j,X}}\ci d(\theta) = \ba{X^{d(\theta)}\ci \pi_{j,X}} = \ba{\pi_{i,X}}$. 
This entails the commutative diagram:
$$
\xymatrix@R=2mm@C=16mm{
d(i)\ar[dd]_-{d(\theta)}\ar[dr]^-{\la\ba{\pi_{i,X}}\ra_X}\\
&\Gamma\\
d(j)\ar[ur]_-{\la\ba{\pi_{j,X}}\ra_X}
}
$$
in view of the fact that $\la\ba{\pi_{j,X}}\ra_X\ci d(\theta) =
\la\ba{\pi_{j,X}}\ci d(\theta)\ra_X$. 

Now that we have constructed our cocone, we have to prove that it is weakly initial
among all cocones on $d$. Let $(\delta_i:d(i)\rar D)_i$ be a cocone on $d$ with
vertex $D$. We have $\delta_j\ci d(\theta) = \delta_i$ for all $\theta:i\rar j$.
We have to construct an arrow $\psi:\Gamma\rar D$, such that $\psi\ci\gamma_i =
\delta_i$ for all $i$. 

We have the arrow $\pi_D:\Gamma\rar D^{\Lim D^d}$. We also  have the element $[\delta]$
in $\Lim D^d$. 
Now, we can define $\psi$ as the following composition:
$$
\xymatrix{
\Gamma\ar[r]^-{\pi_D}
  &D^{\Lim D^d}\ar[r]^-{\ev_{[\delta]}}
     &D
}
$$
We have to check that $\delta_i = \psi\ci\gamma_i$ (for all $i$):
$$
\begin{array}{rcl}
\psi\ci\gamma_i &=& \ev_{[\delta]}\ci\pi_D\ci\la\ba{\pi_{i,X}}\ra_X\\
                &=& \ev_{[\delta]}\ci\ba{\pi_{i,D}}\\
                &=& \delta_i
\end{array}
$$
This completes the proof of the theorem.~\qed

\end{document}